\documentclass[12pt,epsfig]{article}
\usepackage[frenchb]{babel}
\usepackage[T1]{fontenc}
\usepackage[latin1]{inputenc}
\usepackage{amssymb,amsbsy,amsmath,amsfonts,amssymb,amscd}

\begin{document}
\begin{center}
{\large\bf Revêtements hyperelliptiques $d$-osculateurs
 et solitons elliptiques de la hiérarchie $KdV$}\\[5mm]
Armando \textsc{Treibich}\\
Faculté J.Perrin, rue J.Souvraz, SP18, 62300 Lens, France\\
Laboratoire de Mathématique de Lens, Université d'Artois\\
e-mail addresses: treibich@euler.univ-artois.fr ; treibich@cmat.edu.uy
\end{center}

\vspace*{5mm} \noindent\textbf{Résumé -} Soit $d$ un entier positif,
$ \mathbb K$ un corps algébriquement clos de caractéristique 0 et $ X$ une courbe
elliptique définie sur $\mathbb K$. On étudie les courbes
hyperelliptiques munies d'une projection sur $ X$, telles que
l'image naturelle de $X$ dans la jacobienne de la courbe, oscule à
l'ordre $d$ au plongement de celle-ci, en un point de Weierstrass.
On construit des familles ($d-1$)-dimensionnelles de telles courbes,
de genre $g$ arbitrairement grand, obtenant, en particulier, des
familles $(g+d-1)$-dimensionnelles de solutions de la hiérarchie
$KdV$, doublement périodiques par rapport à la $d$-ième variable.

\vspace*{2mm}

\noindent\textbf{Abstract -} Let $d$ be a positive integer, $\mathbb
K$ an algebraically closed field of characteristic 0 and $ X$ an elliptic curve defined
over $\mathbb K$ . We study the hyperelliptic curves equipped with a projection
over $ X$, such that the natural image of $ X$ in the Jacobian of
the curve osculates to order $d$ to the embedding of the curve, at a
Weierstrass point. We construct ($d-1$)-dimensional families of such
curves, of arbitrary big genus $g$, obtaining, in particular,
$(g+d-1)$-dimensional families of solutions of the $KdV$
hierarchy, doubly periodic with respect to the $d$-th variable.\\

\noindent\textbf{Abridged english version}

Let $\mathbb K$ be an algebraically closed field of characteristic 0 and let $\mathbb
P^1$ denote the projective line over $\mathbb K$. By a curve we will
mean hereafter a complete integral curve over $\mathbb K$, of
arithmetic genus $ g>0$. For any curve $\Gamma$, let $\Gamma ^o $
and $\textrm{Jac}\ \Gamma$ denote, respectively, the open subset of
smooth points of $\Gamma$ and its generalized Jacobian. Recall that
for any smooth point $p\in\Gamma ^o$, the Abel morphism, $
A_p:\Gamma ^o \to \textrm{Jac}\ \Gamma$, $ p'\mapsto O_\Gamma
(p'-p)$, is an embedding and $A_p(\Gamma ^o)$ generates the whole
Jacobian. Moreover, the flag of hyperosculating spaces to
$A_p(\Gamma ^o)$ at $ A_p(p)= 0$  can be constructed as follows. For
any marked curve $(\Gamma , p)$ as above, and any positive integer
$j$, let us consider the exact sequence of $O_\Gamma $-modules $0\to
O _\Gamma \to O_\Gamma (jp)\to O_{jp} (jp)\to 0$, as well as the
corresponding long exact cohomology sequence :

$$0\to H^o(\Gamma,O_\Gamma )\to H^o(\Gamma,O_\Gamma (jp))\to H^o(\Gamma,O_{jp}(jp))\stackrel{\delta }\to H¹(\Gamma,O_\Gamma )\to \ldots ,$$
where $\delta : H^o(\Gamma,O_{jp} (jp)) \to H^1(\Gamma,O_\Gamma )$
is the canonical cobord morphism and $H^1(\Gamma,O_\Gamma )$ is
canonically identified to the tangent space to $\textrm {Jac}\
\Gamma $ at $ 0$. According to the Weierstrass gap Theorem, for any
$d=1,\ldots ,g$ ($g$ being the arithmetic genus of $\Gamma $), there
exists $0 <j<2g$ such that $\delta (H^o(\Gamma,O_{jp}(jp)))$ is a
$d$-dimensional subspace, denoted hereafter by  $W_d$. For a generic
point $p$ of $\Gamma $ we have $W_d= \delta (H^o(\Gamma,O_{dp}(dp)))$
(i.e.: $j=d)$, but if $\Gamma $ is hyperelliptic and $p$ is a
Weierstrass point, then we must choose $j=2d-1$.

\noindent In any case, the filtration $\{0\} \subsetneq  W _1\ldots
\subsetneq W _g = H^1(\Gamma,O_ \Gamma)$ is the flag of
hyperosculating spaces to $A_p(\Gamma )$ at $ 0 $.  \\

\noindent Definition 1 (cf.[7])\\
 \textit{Let $\pi:(\Gamma
,p)\to(X,q)$ be a finite marked morphism $(\pi(p)=q)$ and let
$\iota_\pi:X\to \textrm{Jac}\  \Gamma $ denote the canonical (group
homo-)morphism $q'\mapsto A_p(\pi^*(q'-q))$. We will say that $\pi$
is a $d$-osculating cover iff the tangent to $\iota_\pi(X)$ at $ 0 $
is contained in $W _d $ but not in $W _{d-1}$. If, moreover,
$\Gamma$ is a degree-$2$ cover of $\mathbb P^1$ ramified at $p$, we
will call $\pi$ a hyperelliptic $d$-osculating cover.}
\vspace*{2mm}

\noindent{Remark 2 } \\Besides its intrinsic geometric interest, the
above definitions are motivated, when $ \mathbb K = \mathbb C$, by
the following result : any hyperelliptic $d$-osculating cover of
genus $g$ gives rise to a $g$-dimensional family of exact solutions
of the Korteweg-de Vries equation, doubly periodic with respect to
 the $d$-th $KdV$ flow (cf.[3],[4],[6],[5]).\\

\noindent Definition 3 (cf.[8]) \begin {enumerate}
 \item \textit{ Let E
denote the unique indecomposable rank-$2$, degree-$0$ vector bundle
over X and $ S  : =  \mathbb P(E)$ the corresponding projective
bundle. The natural projection $\pi_s : S \to X $ is then a ruled
surface over X, characterized up to an isomorphism, by the existence
of a unique section, $C_o\subset S$, of self-intersection $C_o.C_o =
0$.}

\item  \textit{The canonical symmetry of $(X,q),
[-1]: X \to X$, fixes its origin,
$\omega _o : = q$, as well as the three other half-periods $
\{\omega _j,j = 1,2,3\}$, and lifts to an involution $\tau : S \to
S$ $ (i.e.:\pi_{s}\circ \tau = [-1]\circ \pi_s)$ having two fixed
points over each $\omega _i$ $(i=0,\ldots,3)$ : one in $ C_o$, denoted by
$s_i$, and the other one denoted by $r_i$.}

\item \textit{ Let $e$ : $ S^\bot \to S $ denote the blow-up of $S$
at $\{s_i,r_i,i=0,\ldots,3\}$, the eight fixed points of $ \tau$, and
$\tau ^\bot : \ S^\bot \to S^\bot$ its lift to an involution fixing
the corresponding exceptional divisors
${s_i^\bot,r_i^\bot,i=0,\ldots,3}$. Then, the quotient $\tilde S :
=S^\bot/\tau^\bot$ is a smooth rational surface and the canonical
degree-$2$ projection $ \psi:S^\bot\to \tilde S $ is ramified along
${s_i^\bot,r_i^\bot,i=0,\ldots,3}$.}

\end{enumerate}

- The following results, already known for $d = 1$ (cf.[8]) and $
d = 2$ (cf.[1]), can be proven within the same framework for any $d>2.$\\

\noindent Proposition 4 (cf.[7])\\
 \noindent \textit{Let $\pi : (\Gamma ,p)\to (X,q)$  be a hyperelliptic
$d$-osculating cover of degree  $n$  and let  $\rho$  denote its ramification index
at $p$. Then  $\rho$  is odd, bounded by  $2d-1$, and there exists a
unique morphism   $ \iota^\bot:\Gamma \to  S^\bot $  such that}:
\begin {enumerate}
\item \textit{the surjective morphism  $ \iota^\bot : \Gamma \to \iota ^\bot(\Gamma
)$ factors  $\pi = \pi_{s} \circ e \circ \iota^\bot$  and its degree
divides $ 2d-1$};

\item \textit{the curve $\iota^\bot(\Gamma ) $ is $\tau ^\bot$-invariant
and projects into $ \tilde S $, with degree $2$ over the rational
curve  $\tilde \Gamma \ :=  \varphi (\iota^\bot(\Gamma ))$};

\item \textit{the direct image divisor $ (e \circ \iota^\bot)_*(\Gamma ) $ is linearly
equivalent to $nC_o +(2d-1)S_o$ and only
intersects $C_o$ at $(e \circ \iota^\bot)(p) = C_o\cap S_o$};

\item \textit{ for any $i = 0,\ldots,3$ let  $\gamma _i$  denote the intersection multiplicity number
between the direct image divisor  $ (\iota^\bot)_*(\Gamma )$  and  $
r_i^\bot$. Then  $(\iota^\bot)_*(\Gamma )$  is linearly equivalent
to  $e^*(nC_o+(2d-1)S_o)-\rho s_o^\bot-\Sigma _ i\gamma_i r_i^
\bot$.}

\end{enumerate}
\noindent Definition 5

\noindent \textit{ Let  $\gamma  :=(\gamma _i)\in \mathbb N^4, \gamma _i := (\iota^\bot)_* (\Gamma).r_i^\bot (i=0,\ldots,3)$ be the intersection multiplicity vector canonically associated to the hyperelliptic
$d$-osculating cover $\pi : (\Gamma ,p)\to (X,q)$.
 We will call $\gamma$ the type of $\pi$ and denote hereafter by $\gamma ^{(1)}$ and
$\gamma^{(2)} $ the sums $\gamma ^{(1)} := \Sigma_i \gamma_i$ and
$\gamma ^{(2)} :=  \Sigma_i\gamma_i ^2$, respectively.}\\

\noindent Theorem 6 \\
\noindent \textit{Let $\pi:(\Gamma ,p)\to (X,q)$ be a hyperelliptic
$d$-osculating cover, of degree $n$, type $ \gamma $, ramification
index  $\rho$  at  $p$  and arithmetic genus $g$. Then :}

\begin{enumerate}

\item \textit{$\gamma _o + 1 \equiv \gamma _1\equiv \gamma _2 \equiv \gamma
_3 \equiv \ n \ ($mod$.2)$};

\item \textit{$2g + 1\leq \gamma ^{(1)}$  and  $\gamma^{(2)} \leq (2d-1)(2n-2)+4- \rho^2$}.

\end{enumerate}

\noindent  - Let us fix hereafter  $k\in\{0,1,2,3\}$
and choose $\varepsilon = (\varepsilon_i) \in \mathbb Z^4$  such that, either
$\vert\varepsilon_i\vert = (d-1)(1-\delta
_{i,k})$, for any  $i=0,\ldots ,3$, or $\vert\varepsilon_i\vert = [d/2] -(-1)^d\delta
_{i,k}$,  for any  $i=0,\ldots ,3$. We prove the \\

\noindent Theorem 7

\noindent \textit{For any  $\mu\in \mathbb N^4$  such that
$\mu_o+1\equiv \mu_j($mod.$2) (j=1,2,3)$  and  $ \gamma : =
(2d-1)\mu  +  2\varepsilon \in \mathbb N^4$, there exists a
$(d-1)$-dimensional family of smooth hyperelliptic $d$-osculating covers of
genus  $g : = 1/2\{(2d-1)\mu^{(1)}+2\varepsilon ^{(1)}-1\}$, degree $n: = \frac{1}{2} \bigl \{(2d-1)\mu^{(2)}+4\Sigma
_i\mu_i \varepsilon_i+6d-7 \bigr \}$ and type $\gamma $ .} \\

 \noindent\textbf{Version française}

\noindent 1.1 Soit $\mathbb K$ un corps algébriquement clos de caractéristique 0, $(X,q)$
une courbe elliptique définie sur $\mathbb K$, pointée en son
origine, et notons $\mathbb  P^1$  la droite projective sur $\mathbb
K$. Sauf mention contraire, toutes les courbes considérées par la
suite seront supposées définies sur $\mathbb K$, complètes, intègres
et de genre arithmétique positif. Etant donnée une telle courbe
$\Gamma $, munie du choix d'un point lisse $ p \in \Gamma $, on
désignera par $A_p :\Gamma \to  \textrm{Jac} \Gamma $ l'application
(rationnelle) d'Abel.  \\

\noindent Proposition 1.2 (cf.[7]\S1.6)

\textit{Soit $\Gamma $ une courbe hyperelliptique de genre positif $g$, $p$ un point
lisse de Weierstrass de $\Gamma $ et considérons, quel que soit le
nombre impair $j=2d-1$ < $2g$, la suite exacte de $O_\Gamma$-modules $
0\to O_\Gamma \to O_\Gamma (jp)\to O_{jp}(jp) \to 0$, ainsi que sa
suite exacte
longue de cohomologie}\\
 $0\to H^o(\Gamma,O_\Gamma )\to H^o(\Gamma,O_\Gamma (jp))\to H^o(\Gamma,O_{jp}(jp))
 \stackrel {\delta }{\to }  H^1(\Gamma,O_\Gamma)\to \ldots $,\\

\noindent\textit{ où $\delta $ est l'application cobord. Alors $W
_d=\delta(H^o(O_{jp}(jp)))$  est le $d$-ième sous-espace osculateur
à
$A_p(\Gamma )\ en\ 0.$}\\

\noindent Définition 1.3

\textit{Soit $\pi:(\Gamma ,p)\to (X,q)$ un revêtement ramifié de la
courbe elliptique $(X,q), p$ un point lisse de $\Gamma $ tel que
$\pi(p)=\ q \ $ et notons $\iota _{\pi}: (X,q)\to (Jac \Gamma ,0)$
l'homomorphisme de groupes algébriques, $q'\mapsto
A_p(\pi^*(q'-q))$. Nous dirons que $\pi$ est un revêtement
hyperelliptique $d$-osculateur $(d > 0)$ si et seulement s'il
satisfait les propriétés ci-après} :
\begin{enumerate}
 \item \textit{$\Gamma $
est un revêtement double de $\mathbb P¹$, ramifié en $p$ };

\item \textit{la tangente à $\iota_\pi(X)\ en \ 0$ est contenue dans $W
_d$ mais pas dans $W _{d-1}$}.
\end{enumerate}

\noindent Définition 1.4 (cf.[8])

\begin{enumerate}
 \item \textit{ Soit $\pi_s : S\to X $ l'unique
surface réglée au dessus de X, ayant une seule section, $C_o\subset
S$, d'auto-intersection nulle. La symétrie canonique de $(X,q),
[-1]: X \to X$, fixe l'origine $\omega _o :\ = q$ ainsi que les
trois autres demi-périodes $\{\omega _j, j=1,2,3\}$, et se remonte
en une involution de $S$, notée $\tau : S\to S \ (\pi_s \circ \tau \
= \ [-1]\circ \pi_s)$, ayant deux points fixes au dessus de chaque
$\omega _i(i=0,\ldots ,3)$: un sur $C_o$, noté $s_i$, et l'autre noté
$r_i$.}

\item \textit{Soit d'autre part $e :S^\bot\to S $ l'éclatement des huit
points  $ \{ s_i,r_i \} \subset S $ et notons $\{ s_i^\bot,
r_i^\bot \} $ les diviseurs exceptionnels correspondants.
Alors $\tau : S \to S$ se remonte à son tour en une involution $ \
\tau^ \bot : S^\bot \to S^\bot (e \circ \tau ^\bot = \tau \circ
e)$}.
 \item \textit{ Il s'en suit que la surface
quotient $\tilde S :\ = S^\bot/\tau ^\bot$ est lisse et que le
revêtement double associé,  $
\varphi : S^\bot\to \tilde S$,  est ramifié le long des courbes $\{s_i^\bot,r_i^\bot\}$.}\\

\end{enumerate}

\noindent - Les projections $\pi_s^\bot := \pi_s\circ e :S^\bot \to
X$ et $\varphi : S^\bot \to \tilde S$, donnent un cadre universel
pour les revêtements hyperelliptiques $d$-osculateurs ($d$ > $0$) et permettent de démontrer la proposition 1.6 et le Théorème 1.9 ci-dessous, de façon analogue aux cas
particuliers $d = 1$, $d = 2$ et $\mathbb K = \mathbb C$ (cf.[8],[1]).\\

\noindent Proposition 1.5 (cf.[8],[1],[2])

 \noindent \textit{ Soit $\pi : (\Gamma ,p)\to (X,q)$  un revêtement
hyperelliptique $d$-osculateur de degré n et notons $\rho$ son indice de
ramification en $p$. Alors $\rho$ est un nombre impair majoré par $2d-1$
et il existe un unique morphisme  $\iota ^\bot : \Gamma \to S^ \bot
$ tel que }:
\begin{enumerate}
\item \textit{la courbe image $\iota ^\bot(\Gamma ) \subset S ^\bot $ est
$\iota ^\bot$- invariante et sa projection, $\tilde \Gamma : \ =
\varphi (\iota ^\bot(\Gamma ))$, est une courbe rationnelle
irréductible de $\tilde S$};

\item \textit{$\pi$ se factorise via $\iota^\bot$, $\pi = \pi_s\circ e \circ
\iota^\bot$, et $deg(\iota^\bot : \Gamma \to \Gamma ^\bot)$ divise
$2d-1$};

\item \textit{le diviseur image directe $ (e \circ \iota^\bot)_*(\Gamma ) $ est linéairement équivalent à $nC_o +(2d-1)S_o$ et n'intersecte
  $C_o$ qu'au point $(e \circ \iota^\bot)(p) = C_o\cap S_o$ };

\item \textit{Le diviseur  $(\iota^\bot)_*(\Gamma )$  est linéairement équivalent à
$e ^
*(nC_o+(2d-1)S_o)-\rho s_o ^\bot - \Sigma _i\gamma _ir_i^\bot$, où, pour tout $i=0,\ldots ,3$,  $\gamma _i$ $ := \iota ^\bot_*(\Gamma ).r_i^\bot$ }.
 
\end{enumerate}
 \noindent Définition 1.6

\noindent \textit{Le vecteur $\gamma =(\gamma _i)\in \mathbb N^4,$ $
\gamma _i :\ =$ $ (\iota^\bot)_*(\Gamma ).r_i^\bot$,  canoniquement
associé au revêtement hyperelliptique $d$-osculateur $\pi$, sera
appelé le type de $\pi$. Nous désignerons dorénavant, par  $\gamma
^{(1)}$  et $\gamma
^{(2)}$  les sommes  $\Sigma _i\gamma _i$  et  $\Sigma _i\gamma_i^2$, respectivement.}\\

\noindent Theorème 1.7

\noindent \textit{Soit $\pi:(\Gamma ,p)\to (X,q)$ un revêtement
hyperelliptique $d$-osculateur, de degré $n$, type $\gamma $,
genre arithmétique $g$ et notons  $\rho$  son indice de ramification en p. Alors }:

\begin{enumerate}
\item \textit{$\gamma_o+ 1\equiv\gamma _1\equiv\gamma _2\equiv\gamma _3\equiv n($mod.$2)$};

\item \textit{$2g+1\leqslant\gamma ^{(1)}$ et $\gamma ^{(2)}\leqslant (2d-1)(2n-2)+4-\rho^2$}.

\end{enumerate}
\noindent - Au moyen des critères d'existence et d'irréductibilité ci-après, nous construisons finalement, pour tout $n$ et pour $g$
arbitrairement grand, une famille de dimension $d-1$ de
tels revêtements, de degré $n$ et genre $g$.\\

\noindent Proposition 1.8 (cf.[8]\S6.2)

\noindent \textit{Soit $k\in \{0,1,2,3\}$ et $\alpha  = (\alpha
_i)\in \mathbb N^4$  tel que  $\alpha _k+1 \equiv\alpha
_j$(mod.$2$), quel que soit $j\neq k$. Alors, il existe une unique
courbe irréductible tracée dans $S^\bot$, notée ci-après  $Z^\bot
_\alpha$, qui soit $\tau ^\bot$-invariante et linéairement
équivalente à
 $e^*(nC_o+S_k)-s_k^\bot-\Sigma _i\alpha _ir_i^\bot$,  où  $n$ $ :=
1/2(-1 +\Sigma_i\alpha _i^2)$.}\\

\noindent Proposition 1.9 (cf.[7]\S3.4)

\noindent\textit{ Soit $\Gamma $ un diviseur effectif de la surface
$S$, lisse
en $s_o$ et tel que $\Gamma \cap C_o=s_o$. Alors $\Gamma $ est une courbe irréductible.}\\

\noindent Proposition 1.10

\noindent \textit{Soit $C_o^\bot$ le transformé strict de $C_o$ dans
$S^\bot$ et $\Gamma ^\bot$ un diviseur effectif de $S^\bot$
satisfaisant les propriétés suivantes }:
\begin{enumerate}

\item \textit{la courbe $\Gamma ^\bot$  intersecte $C_o^\bot$  uniquement au point $p_o^\bot :\ = C_o^\bot \cap s_o^\bot$ et son support ne contient aucune des courbes dans  $\{C^\bot_o, s_i^\bot, r_i^\bot, i= 0,\ldots ,3\}$};

\item \textit{quel que soit $i=0,\ldots,3, $   deg$(\Gamma ^\bot.s_i^\bot) = \delta _{i,0}$} . 
\end{enumerate}
\textit{ Alors $\Gamma ^\bot$ est une courbe irréductible }.\\
\noindent Preuve

La propriété 1.10.1. nous assure que $\Gamma ^\bot$ est le
transformé strict de  $\Gamma \ :\ =\ e _*(\Gamma ^\bot)$, son image
directe par  $e:S^\bot \to S$, et que celle-ci ne contient pas
$C_o$. On vérifie également, grâce aux autres propriétés, que
$\Gamma $ est lisse en $s_o$ et que  $\Gamma \cap C_o=s_o$. Il s'en
suit, d'après [7]\S3.4, que $\Gamma $ est une courbe irréductible,
de même que $\Gamma^\bot$, son transformé strict.$  \blacksquare $\\

\noindent Theorème 1.11

\noindent \textit{Fixons  $d\in \mathbb N ,\ \mu = (\mu_i)\in
\mathbb N ^ 4$ tel que  $\mu _o+1\equiv \mu_1\equiv \mu_2\equiv
\mu_3$(mod.$2$) et choisissons  $\varepsilon \in \mathbb Z^4$  égal,
à moins d'une permutation ou d'un changement des signes des
coefficients, soit à $\varepsilon =(0,\ 2d-2,\ 2d-2,\ 2d-2)$, soit à
$\varepsilon = (d-2,d,d,d)$  si  $d$ est pair ou à  $
\varepsilon=(d+1,\ d-1,\ d-1,\ d-1)$  si $d$ est impair. Notons
$\gamma :\ =(2d-1)\mu + \varepsilon$  et soit  $n\in\mathbb N^*$
l'unique naturel positif tel que  $\gamma ^ {(2)} \ := \ \Sigma
_i\gamma _i^2 \ = \ (2d-1)(2n-2)+3$.  Alors, le système linéaire
$\mid e^
*(nC_o+(2d-1)S_o)-s_o^\bot- \Sigma _i \gamma _i r_i^\bot  \mid $
contient un sous-espace de dimension $(d-1)$, dont l'élément
générique est une courbe irréductible, $\tau ^\bot$-invariante et
lisse au point $p_o^\bot :\ = C_o^\bot \cap s_o^\bot$, telle que sa
projection dans $\tilde S$ est isomorphe à
$\mathbb P^1$.}\\

\noindent Corollaire 1.12

\noindent \textit{Soient  $\gamma \in \mathbb N^4$  et  $n\in\mathbb
N$ (où
 $\gamma ^{(2)}= (2d-1)(2n-2)+3))$ comme ci-dessus. Il existe alors
une famille $(d-1)$-dimensionnelle de revêtements hyperelliptiques
$d$-osculateurs, de degré $n$, type $\gamma $ et non-singuliers de genre  
$ g $ $: =1/2(-1+\gamma ^{(1)}).$}\\

\noindent Preuve du Théorème.

\noindent Pour des raisons d'espace nous
allons construire uniquement la famille associée à  $\gamma =(2d$-$1)\mu + (0,\
2d$-$2, \ 2d$-$2, \ 2d$-$2)$, auquel cas le degré $n$ correspondant est tel que  $2n=(2d$-$1) \mu ^{(2)} + 4(d$-$1)(\mu_1 +
\mu_2 + \mu_3)+6d-7\ $.

\noindent Soient  $\mu (1) :\  =
 \mu + (0,0,1,1)$,  $\mu(2):\ = \mu + (0,1,0,1)$,  $\mu(3): \ =
\mu + (0,1,1,0)$,  et notons, d'après 1.11 ci-dessus,  $ Z^\bot_{(k)}  (k=1,2,3)$  l'unique courbe
$\tau^\bot$-invariante de $S^\bot $, telle que   $Z^\bot_{(k)} \sim
e^*(m(k)C_o+S_k)-s_k^\bot - \Sigma _i\mu(k)_ir_i^\bot$,  où   $
2m(k)+ 1= \Sigma_i\mu(k)_i^2$.

\noindent Soient d'autre part  $\mu^+:=\mu+(1,1,1,1)$, $ \mu' :=
\mu+(0,2,1,1)$,   et notons  $Z_+^\bot ,
{Z'}^\bot$,   les uniques courbes  $\tau
^\bot$-invariantes de $ S ^\bot $, linéairement équivalentes à :
\begin{enumerate}

\item $Z^\bot_+ \sim e^*(m^+C_o+S_o)-s_o^\bot-\Sigma _i\mu_i^+r_i
^\bot $, où   $2m^+ + 1 = \mu^{+(2)}$ ;

\item $Z'^\bot \sim e^*(m'C_o+S_o)-s_1^\bot-\Sigma _i\mu_i'
r_i^\bot$, où  $2m'+1=\mu'^{(2)}$.

\end{enumerate}
\noindent De même, si  $\mu_o\neq 0$  on note  $\mu^- $ $:  = $ $\mu
+ (-1,1,1,1)$  et  $Z_-^\bot$  l'unique courbe  $\tau
^\bot$-invariante de $S^\bot$, telle que  $ Z_-^\bot \sim
e^*(m^-C_o+S_o)-s_o^\bot -\Sigma _i\mu_i^-r_i^\bot$,  où   $2m^-+1=
\mu^{-(2)}$.
\noindent Par contre, si  $\mu_o=0$  on notera  $Z_-^\bot :\  = Z_ +
^\bot + 2r_o^\bot$,  de telle sorte que dans les deux cas de figure,
les diviseurs    $D^\bot_o : = Z ^\bot _ +  +  Z^\bot _- + 2s^\bot
_o$   et   $ D^\bot_1 : = Z'^\bot + Z^\bot_{(1)}+ 2s_1^\bot$   soient
linéairement équivalents. 

\noindent Considérons finalement la courbe   $Z^\bot \sim
e^*(mC_o+S_o) - s^\bot_o- \Sigma_i\mu_ir_i^\bot$   où   $2m+1 =
\Sigma_i \mu_i^2$ et remarquons les faits suivants :
\begin{enumerate}
\item tout élément $\tau^\bot$-invariant de $ \Lambda$ est l'image
réciproque par $ \varphi : S^\bot \to \tilde S $, d'un diviseur de
genre arithmétique nul de $ \tilde S$ ;

\item le diviseur    $F^\bot_j : = C^\bot_o + \Sigma _{k\ne 0}(Z^\bot_{(k)}+s^\bot_k) + jD^\bot_o + (d-2-j)D^\bot_1$, pour tout $j=0,\ldots,d-2$, ainsi que  $G^ \bot : = Z^\bot + (d-1)D^\bot_o$, appartiennent à $\Lambda$;

\item les diviseurs $\{F^\bot_j,j=0, \ldots,d-2\}$ ont $C^\bot_o$
comme unique composante irréductible commune, mais $F^\bot_o$ est le
seul à être lisse au point $p^\bot_o :  = C^\bot_o \cap s_o^\bot$.
En particulier ils engendrent un sous-espace $(d-2)$-dimensionnel de
$\Lambda$, dont l'élément générique est
transverse à $s_o^\bot$ en $p^\bot_o$;

\item  le diviseur $G^\bot$ intersecte $C^\bot_o$
uniquement au point $p^\bot_o$.
\end{enumerate}
\noindent Il en résulte que $\{G^\bot, F^\bot _j,j=0,\ldots, d-2\}
$  engendre un sous-espace $(d-1)$-dimensionnel et $\tau^\bot$-invariant de $ \Lambda$, dont
l'élément générique satisfait le critère d'irréductibilité 1.10
et se projette dans $ \tilde S$ sur une
courbe isomorphe à $\mathbb P^1$. $  \blacksquare $\\

\noindent Preuve du Corollaire :

\noindent Soit  $ \Gamma ^\bot $  l'élément générique du sous-espace de  $ \Lambda$
 construit
ci-haut. On sait, d'après le Théorème 1.11, que $ \varphi : \Gamma^\bot \to \tilde \Gamma $  est un revêtement de degré 2, ramifié en $p_o^\bot$, dont l'image est isomorphe à  $\mathbb P^1$. Donc   $\Gamma ^\bot$  est une courbe
hyperelliptique et $p_o^\bot \in \Gamma ^\bot$ est un point de
Weierstrass de   $\Gamma ^\bot$. Il s'en suit que la
projection naturelle de $ (\Gamma ^\bot,p_o^\bot)$ sur $(X,q)$
(restriction de  $\pi_s^\bot :S^\bot \to X $  à  $\Gamma ^\bot)$,
est un revêtement hyperelliptique $d$-osculateur de type $\gamma $,
degré $n$ et genre $g$, tel que  $(2n-2)(2d-1)+3=\gamma ^ {(2)} $
  et   $2g +1=\gamma ^{(1)}.$  $  \blacksquare $


\end{document}